\documentclass[11pt]{article}
\usepackage{stmaryrd}
\usepackage{mathrsfs}
\usepackage[centertags]{amsmath}
\usepackage{amsfonts, dsfont}
\usepackage{amssymb}
\usepackage{amsthm}
\usepackage{newlfont}
\usepackage{bezier,  amsxtra, amsbsy, amsgen, amsopn, amstext}
\usepackage{verbatim}
\usepackage{appendix}

\input xy
\xyoption{all}

\theoremstyle{plain}
\newtheorem{thm}{Theorem}[section]
\newtheorem{cor}[thm]{Corollary}
\newtheorem{lem}[thm]{Lemma}
\newtheorem{prop}[thm]{Proposition}

\theoremstyle{definition}
\newtheorem{defn}[thm]{Definition}

\newtheorem*{remark}{\textbf{Remark}}
\newtheorem*{remarks}{\textbf{Remarks}}
\newtheorem*{ack}{Acknowledgments}

\newcommand{\bd}{\begin{defn}}
\newcommand{\ed}{\end{defn}}
\newcommand{\bl}{\begin{lem}}
\newcommand{\el}{\end{lem}}
\newcommand{\bp}{\begin{prop}}
\newcommand{\ep}{\end{prop}}
\newcommand{\bt}{\begin{thm}}
\newcommand{\et}{\end{thm}}
\newcommand{\bc}{\begin{cor}}
\newcommand{\ec}{\end{cor}}
\newcommand{\br}{\begin{remarks}}
\newcommand{\er}{\end{remarks}}
\newcommand{\bdi}{\begin{diagram}}
\newcommand{\edi}{\end{diagram}}
\newcommand{\beq}{\begin{eqn}}
\newcommand{\eeq}{\end{eqn}}
\newcommand{\ba}{\begin{array}}
\newcommand{\ea}{\end{array}}
\newcommand{\bpf}{\begin{proof}}
\newcommand{\epf}{\end{proof}}

\newcommand{\N}{\mathds{N}}
\newcommand{\Z}{\mathds{Z}}
\newcommand{\Q}{\mathds{Q}}
\newcommand{\Zp}{\mathds{Z}_{p}}
\newcommand{\Qp}{\mathds{Q}_{p}}
\newcommand{\al}{\alpha}
\newcommand{\be}{\beta}

\newcommand{\Ga}{\Gamma}
\newcommand{\ga}{\gamma}
\newcommand{\La}{\Lambda}
\newcommand{\la}{\lambda}
\newcommand{\bb}{\mathbf{b}}

\newcommand{\m}{\mathfrak{m}}
\newcommand{\M}{\mathfrak{M}}

 \DeclareMathOperator{\Sel}{Sel}
\DeclareMathOperator{\Gal}{Gal} \DeclareMathOperator{\Hom}{Hom}
\DeclareMathOperator{\Ext}{Ext} \DeclareMathOperator{\Cl}{Cl}
\DeclareMathOperator{\Tor}{Tor} \DeclareMathOperator{\rank}{rank}
\newcommand{\ot}{\otimes}

\newcommand{\ilim}{\displaystyle \mathop{\varinjlim}\limits}
\newcommand{\plim}{\displaystyle \mathop{\varprojlim}\limits}
\newcommand{\im}{\mathrm{im}\,}

\newcommand{\cts}{\mathrm{cts}}

\newcommand{\cd}{\mathrm{cd}}
\newcommand{\cyc}{\mathrm{cyc}}

\newcommand{\lra}{\longrightarrow}
\newcommand{\tha}{\twoheadrightarrow}

\newcommand{\sbs}{\subseteq}

\newcommand{\ps}[1]{\llbracket #1 \rrbracket}

\begin{document}

\title{Notes on the fine Selmer groups}
\author{Meng Fai Lim\footnote{Department of Mathematics, University of
Toronto, 40 St.\ George St., Toronto, ON, M5S 2E4  Canada, e-mail:
mengfai.lim@utoronto.ca} \footnote{Present Address: School of
Mathematics and Statistics, Central China Normal University, 152
Luoyu Road, Wuhan, Hubei, P.R.China 430079, e-mail:
limmf@mail.ccnu.edu.cn}}
\date{}
\maketitle

\begin{abstract} \footnotesize
\noindent  In this paper, we study the fine Selmer groups attached
to a Galois module defined over a commutative complete Noetherian
ring with finite residue field of characteristic $p$. Namely, we are
interested in its properties upon taking residual representation and
within field extensions. In particular, we will show that the
variation of the fine Selmer group in a cyclotomic $\Zp$-extension
is intimately related to the variation of the class groups in the
cyclotomic tower.


\end{abstract}

\normalsize
\section{Introduction}

Let $F$ be a finite extension of $\Q$, $p$ a prime and $F^{\cyc}$
the cyclotomic $\Zp$-extension of $F$. Denote $K(F^{\cyc})$ to be
the maximal unramified pro-$p$ extension of $F^{\cyc}$ at which
every prime of $F^{\cyc}$ above $p$ splits completely. Set $\Ga =
\Gal(F^{\cyc}/F)$. Iwaswa has proven that
$\Gal(K(F^{\cyc})/F^{\cyc})$ is a finitely generated torsion
$\Zp\llbracket\Ga\rrbracket$-module. He further conjectured that it
is in fact a finitely generated $\Zp$-module (see \cite{Iw, Iw2}).
Conjectures parallel to this conjecture of Iwasawa have been
formulated for the fine Selmer groups attached to elliptic curves by
Coates and Sujatha \cite[Conjecture A]{CS}, and for the fine Selmer
groups attached to modular forms and Hida families by Sujatha and
Jha \cite[Conjecture A, Conjecture 1]{JhS, Jh}. Following their
footsteps, we make an analogous conjecture for the fine Selmer
groups attached to a Galois module defined over a commutative
complete Noetherian ring with finite residue field of characteristic
$p$. We then show that this generalized conjecture turns out to be a
consequence of the original conjecture of Iwasawa (see Theorem
\ref{main} and \ref{main prop}). We remark that such an implication
has been established for elliptic curves by Coates and Sujatha (see
\cite[Theorem 3.4, Corollary 3.5]{CS}) and our approach is mainly
inspired by their work. The proof adopted by them made used of
descent technique and relied on the observation that the extension
carved out by all the $p$-power torsion points of the elliptic curve
in question is a $p$-adic Lie extension. It is not difficult to
extend their argument to a $p$-adic representation defined over a
ring of integers of some finite extension of $\Qp$. However, since
the Galois modules we considered are general, the extensions carved
out by these modules need not be $p$-adic Lie extensions, and so
their argument does not carry over. Therefore, our proof takes a
different route, which we will explain in a while.


In the second part of the article, we will investigate the question
of pseudo-nullity of the fine Selmer groups over an admissible
$p$-adic Lie extension of dimension strictly greater than one. A
somewhat related question in this direction was first considered by
Greenberg \cite{G}, where he conjectured that the Galois group of
the maximal abelian unramified pro-p-extension of the compositum of
all $\Zp$-extensions $\tilde{F}$ of $F$ is a pseudo-null
$\Zp\ps{\Gal(\tilde{F}/F)}$-module. Now if $T$ is the Tate module of
all the $p$-power roots of unity, the dual fine Selmer group is
precisely $\Gal(K(F_{\infty})/F_{\infty})$, where we denote
$K(F_{\infty})$ to be the maximal unramified pro-$p$ extension of
$F_{\infty}$ at which every prime of $F_{\infty}$ above $p$ splits
completely. Hachimori and Sharifi \cite{HS} has constructed several
classes of admissible $p$-adic Lie extensions $F_{\infty}$ of $F$ of
dimension $>1$ such that $\Gal(K(F_{\infty})/F_{\infty})$ is not
pseudo-null as a $\Zp\ps{\Gal(F_{\infty}/F)}$-module. On the other
hand, Coates and Sujatha have conjectured that pseudo-nullity for
the fine Selmer groups of elliptic curves should hold for any
admissible $p$-adic extensions of dimension $>1$ (see
\cite[Conjecture B]{CS}). This was further extrapolated and
formulated for fine Selmer groups of modular forms and Hida families
by Jha \cite[Conjecture B, Conjecture 2]{Jh}.  In this paper, we
will establish two results pertaining to the pseudo-nullity property
of the fine Selmer groups of a Galois module. The first result
(Theorem \ref{pseudo-null main}) shows that the pseudo-nullity
property can be lifted, namely, if the fine Selmer group of the
residual representation of a given Galois representation is
pseudo-null, then so is the fine Selmer group of the Galois
representation itself. Such a result has also been proved in
\cite[Theorem 10]{Jh}\footnote{ Jha's theorem also deals with the
converse direction which we will not address in this paper.}, and
here we give a slightly different proof. The second result (Theorem
\ref{pseudo-null descent}) is concerned on the descent property of
the pseudo-nullity of the fine Selmer groups, where we show that if
$F_{\infty}'\sbs F_{\infty}$ are two admissible $p$-adic Lie
extensions of dimension $\geq 2$ with $\Gal(F_{\infty}/F_{\infty}')$
being a solvable uniform pro-$p$ group, then the fine Selmer group
over $F_{\infty}$ is pseudo-null whenever the fine Selmer group over
$F_{\infty}'$ is pseudo-null. We note that in this result, the
extensions $F_{\infty}$ and $F_{\infty}'$ themselves need not be
solvable extensions of $F$.

We briefly summarize the the approach towards the investigation of
the fine Selmer group in this paper. Namely, we observe that the
structure of the dual fine Selmer group is intimately linked with a
certain inverse limit of second cohomology groups. This latter group
turns out to behave well upon taking residual (see Lemma
\ref{residual H2}) and descent (cf. Lemma \ref{descent H2}), which
is key to our examination of the relationship between the said fine
Selmer groups.

We now give an outline of the paper. In Section \ref{Fine Selmer
groups}, we introduce the fine Selmer groups and discuss some of
their basic properties. In Section \ref{Cyclotomic Zp extension}, we
will study the variation of the fine Selmer groups over the
cyclotomic $\Zp$-extension as mentioned in the first paragraph of
the introduction. In Section \ref{Ranks of Iwasawa modules}, we
collect some results on the ranks of modules over a completed group
algebra which will be applied in Section \ref{pseudo-nullity
section} for the discussion of the pseudo-nullity of the fine Selmer
groups. In Section \ref{example}, we discuss some numerical examples
coming from class groups, and an elliptic curve and its associated
Hida deformation to illustrate the results of the paper. In Section
\ref{Torsion section}, we will make some complementary remark on the
torsionness of the fine Selmer group. Finally, we will provide
proofs to certain results on the structure of the completed group
algebra in the Appendix.

\begin{ack}
    The author would like to thank Romyar Sharifi for his comments and suggestions.
    The author would also like to thank Ramdorai Sujatha for her interest and comments.
    This work was written up when the author is a Postdoctoral fellow at the GANITA Lab
    at the University of Toronto. He would like to acknowledge the
    hospitality and conducive working conditions provided by the GANITA
    Lab and the University of Toronto. He also thank the anonymous referee for pointing out some mistakes and for giving various
    helpful comments that have improved the exposition of the paper.
    \end{ack}

\section{Fine Selmer group} \label{Fine Selmer groups}

Throughout the paper, $p$ will denote a prime number. If $M$ is a
pro-$p$ group or a discrete $p$-primary group, we denote the
Pontryagin dual of $M$ by $M^{\vee}= \Hom_{\cts}(M, \Qp/\Zp)$. Let
$F$ be a number field, i.e., a finite extension of $\Q$. Let $S$
denote a finite set of primes of $F$ containing the primes above $p$
and the infinite primes, which we shall fix once and for all. Let
$F_S$ denote the maximal algebraic extension of $F$ unramified
outside $S$. For any algebraic (possibly infinite) extension
$\mathcal{L}$ of $F$ contained in $F_S$, we write $G_S(\mathcal{L})
= \Gal(F_S/\mathcal{L})$.  Let $R$ be a commutative complete
Noetherian local ring with maximal ideal $\m$ and residue field $k$,
where $k$ is finite of characteristic $p$. Let $T$ denote a finitely
generated $R$-module. One verifies easily that $T\cong \plim_i
T/\m^iT$, and we shall endow $T$ with the $\m$-adic topology. We
assume further that $T$ has a continuous $R$-linear $G_F$-action,
and is unramified outside $S$. This group action on $T$ induces a
continuous group homomorphism
\[ \rho: G_S(F) \lra \mathrm{Aut}_{R}(T). \]
 We will write $W = T^{\vee}(1)$, where $``(1)"$ denotes the Tate
twist. Let $v$ be a prime in $S$. For each finite extension $L$ of
$F$ contained in $F_S$, we define
 \[K_v^i(W/L) = \bigoplus_{w|v}H^i(L_w, W) \quad(i = 0,1),\]
where $w$ runs over the (finite) set of primes of $L$ above $v$. If
$\mathcal{L}$ is an infinite extension of $F$ contained in $F_S$, we
define
\[K_v^i(W/\mathcal{L}) = \ilim_L K_v^i(W/L),\]
where the direct limit is taken over all finite extensions $L$ of
$F$ contained in $\mathcal{L}$ under the restriction maps.

For any algebraic (possibly infinite) extension $\mathcal{L}$ of $F$
contained in $F_S$, the fine Selmer group of $W$ over $\mathcal{L}$
(with respect to $S$) is defined to be
\[ R_S(W/\mathcal{L}) = \ker\Big(H^1(G_S(\mathcal{L}), W)\lra \bigoplus_{v\in S}K_v^1(W/\mathcal{L})
\Big). \] We shall write $Y_S(T/\mathcal{L})$ for the Pontryagin
dual $R_S(W/\mathcal{L})^{\vee}$ of the fine Selmer group. It
follows from the Poitou-Tate sequence that we have the following
exact sequence $$ 0\lra Y_S(T/\mathcal{L}) \lra \plim_{L}H^2(G_S(L),
T)\lra \Big(\bigoplus_{v\in S}K_v^0(W/\mathcal{L})\Big)^{\vee} \lra
W(\mathcal{L})^{\vee}\lra 0,$$ where the inverse limit is taken over
all finite extensions $L$ of $F$ contained in $\mathcal{L}$ under
the corestriction maps and $W(\mathcal{L})$ denote
$W^{\Gal(F_S/\mathcal{L})}$. From now on, we shall write
$H^2_S(\mathcal{L}/F, T) = \plim_{L}H^2(G_S(L), T)$.

\medskip
\noindent Examples: (a) If $T = \Zp(1)$, then we have $W = \Qp/\Zp$.
In this case, one checks easily that $Y_S(\Zp(1)/F) =
\Cl_S(F)[p^{\infty}]$, where $\Cl_S(F)$ is the $S$-class group of
$F$.

\medskip (b) If $E$ is an elliptic curve over $F$, then the fine
Selmer group for $W = E[p^{\infty}]$ has been studied in \cite{CS}
and is related to the classical Selmer group via the following exact
sequence
\[ 0\lra R_S(E[p^{\infty}]/F) \lra \Sel_{p^{\infty}}(E/F) \lra
\bigoplus_{v|p}H^1(F_v, E[p^{\infty}]).\] Note that in this
instance, the fine Selmer group is independent of the set $S$ as
long as it contains all the primes above $p$, the infinite primes
and the primes at which $E$ has bad reduction.

\bigskip
We mention certain basic properties on the structure of the dual
fine Selmer group. Let $F_{\infty}$ be a $p$-adic Lie extension of
$F$ contained in $F_S$. Set $G = \Gal(F_{\infty}/F)$. It is known
that the ring $R\ps{G}$ is Noetherian (cf.\ \cite[Proposition
3.0.1]{LS}, see also Proposition \ref{compact Auslander regular}(a)
in this paper). By \cite[Proposition 4.1.3]{LS}, we have that
$H_S^2(F_{\infty}/F, T)$ is finitely generated over $R\ps{G}$.
Therefore, it follows from the Poitou-Tate sequence that
$Y_S(T/F_{\infty})$ is finitely generated over $R\ps{G}$. We end the
section with two lemmas on the behavior of $H^2_{S}(F_{\infty}/F,
T)$ under residue and descent.

\bl \label{residual H2} Let $x$ be a nonzero and nonunital element
of $R$. Write $\bar{R} = R/xR$ and $\bar{T} = T/xT$. Suppose that
either $p$ is odd or $F$ has no real places. Then one has an
isomorphism
\[ H^2_{S}(F_{\infty}/F, T)/x \cong H^2_{S}(F_{\infty}/F, \bar{T})\]
of $\bar{R}\ps{\Gal(F_{\infty}/F)}$-modules. \el

\bpf Since we are assuming that either $p$ is odd or $F$ has no real
places, we have $H^i_{S}(F_{\infty}/F,-) =0$ for $i\geq 3$. In
particular, $H^2_{S}(F_{\infty}/F,-)$ is right exact. Therefore,
from the exact sequence
\[ T \stackrel{x}{\lra} T \lra \bar{T}\lra 0, \]
we obtain
\[H^2_{S}(F_{\infty}/F, T) \stackrel{x}{\lra} H^2_{S}(F_{\infty}/F, T) \lra H^2_{S}(F_{\infty}/F, \bar{T})\lra 0. \]
The required isomorphism is now immediate.
 \epf

\bl \label{descent H2} Let $F_{\infty}'$ be another $p$-adic
extension of $F$ such that $F\sbs F_{\infty}'\sbs F_{\infty}$
Suppose that either $p$ is odd or $F$ has no real places. Then one
has an isomorphism
\[ H^2_{S}(F_{\infty}/F, T)_{\Gal(F_{\infty}/F_{\infty}')} \cong H^2_{S}(F_{\infty}'/F, T)\]
of $R\ps{\Gal(F_{\infty}'/F)}$-modules. \el

\bpf By \cite[Theorem 3.1.8]{LS}, there is a spectral sequence
 $$ H_i(\Gal(F_{\infty}/F_{\infty}'), H^{-j}_S(F_{\infty}/F,
 T))\Longrightarrow H^{-i-j}_S(F_{\infty}'/F, T). $$
 The required isomorphism follows from reading off the $(0,-2)$-term.
\epf

\begin{remark}
If $T$ is free over $R$, one can deduce Lemmas \ref{residual H2} and
\ref{descent H2} from \cite[Proposition 1.6.5(iii)]{FK}.
\end{remark}

\section{Fine Selmer groups over a cyclotomic $\Zp$-extension} \label{Cyclotomic Zp extension}

In this section, we examine the variation of the fine Selmer groups
over a cyclotomic $\Zp$-extension. We retain the notation and
assumptions from the previous section. \textit{We shall also always
assume that the number field $F$ has no real primes when $p=2$}.
Denote $F^{\mathrm{cyc}}$ to be the cyclotomic $\Zp$-extension of
$F$ and write $\Ga = \Gal(F^{\cyc}/F)$. As mentioned in the previous
section, the dual fine Selmer group $Y_S(T/F^{\cyc})$ is a finitely
generated $R\llbracket\Ga\rrbracket$-module. In fact, we claim that
the following stronger statement should hold.

\medskip \noindent \textbf{Conjecture A.} \textit{For
any number field $F$,\, $Y_S(T/F^{\cyc})$ is a finitely generated
$R$-module.}

\medskip
We mention several well known cases of the above conjecture. For any
extension $\mathcal{L}$ of $F^{\cyc}$ contained in $F_S$, we denote
$K(\mathcal{L})$ to be the maximal unramified pro-$p$ extension of
$\mathcal{L}$ where every prime of $\mathcal{L}$ above $p$ splits
completely. Since $F^{\cyc}\sbs \mathcal{L}$, it follows that every
finite primes of $\mathcal{L}$ splits completely in
$K(\mathcal{L})$. Therefore, in the case when $T = \Zp(1)$, the dual
of the fine Selmer group $Y_S(\Zp(1)/F^{\cyc})$ is precisely
$\Gal(K(F^{\cyc})/F^{\cyc})$. In this context, this is equivalent to
the conjecture made by Iwasawa \cite{Iw, Iw2}. We shall call this
conjecture the Iwasawa $\mu$-invariant conjecture for $F^{\cyc}$.
Currently, the Iwasawa $\mu$-invariant conjecture is only proved in
the case when $F$ is abelian over $\Q$ (see \cite{FW, Sin}).

This conjecture has also been considered for an elliptic curve (see
\cite[Conjecture A]{CS}), for a Galois representation attached to a
cuspidal eigenform which is ordinary at $p$ (see \cite[Conjecture
A]{Jh, JhS}) and a Galois representation attached to a Hida family
(see \cite[Conjecture 1]{Jh}).

\medskip
We now state the main theorem of this section which is a natural
extension of \cite[Corollary 3.5]{CS}.

\bt \label{main} Let $F$ be a number field. Suppose the Iwasawa
$\mu$-invariant conjecture holds for $L^{\cyc}$ for any finite
extension $L$ of $F$. 
  Then $Y_S(T/F^{\cyc})$ is a finitely generated $R$-module. \et

As one will see from the proof (see also Theorem \ref{main prop}),
one only requires that the Iwasawa $\mu$-invariant conjecture holds
for a particular extension $L$ of $F$. We should mention that our
method of proof do not allow us to deduce the finite generation of
$Y_S(T/F^{\cyc})$ from the validity of the Iwasawa $\mu$-invariant
conjecture for $F$, and one has to assume that the Iwasawa
$\mu$-invariant conjecture holds for an extension $L$ of $F$ in
general. In preparation for the proof of Theorem \ref{main}, we
first prove three lemmas. Recall that we write $W(\mathcal{L}) =
W^{\Gal(F_S/\mathcal{L})}$ for any $F\sbs \mathcal{L}\sbs F_S$.
Similarly, for each $v\in S$, we write $W(\mathcal{L}) =
W^{\Gal(\bar{F}_v/\mathcal{L})}$ for any $F_v\sbs \mathcal{L}\sbs
\bar{F}_v$.

\bl \label{fg finite}  Let $F$ be a number field. Suppose that $L$
is a finite extension of $F$ contained in $F_S$. If
$Y_S(T/L^{\cyc})$ is a finitely generated $R$-module, then
$Y_S(T/F^{\cyc})$ is a finitely generated $R$-module. \el

\bpf Let $\triangle = \Gal(L^{\cyc}/F^{cyc})$. Consider the
following commutative diagram with exact rows
\[ \SelectTips{eu}{}
\xymatrix{
  0 \ar[r] & R_S(W/F^{\cyc}) \ar[d]^{\al} \ar[r]^{} & H^1(G_S(F^{\cyc}), W) \ar[d]^{\be}
  \ar[r]^{} & \bigoplus_{v\in S} K_v^1(W/F^{\cyc}) \ar[d]^{\ga} \\
  0 \ar[r] & R_S(W/L^{\cyc})^{\triangle} \ar[r] & H^1(G_S(L^{\cyc}), W)^{\triangle}
  \ar[r] & \big(\bigoplus_{v\in S} K_v^1(W/L^{\cyc})\big)^{\triangle}     }
\]
where the vertical maps are the natural restriction maps. Now $\ker
\al$ is contained in $H^1(\triangle, W(L^{\cyc}))$ which can be
easily seen to be a cofinitely generated $R$-module. The conclusion
is now immediate. \epf

\bl \label{fine selmer and class group} $($Compare with \cite[Lemma
$3.8$]{CS}$)$ Suppose that $F$ contains $\mu_p$ and suppose that $M$
is a finite trivial $G_S(F^{\cyc})$-module which is killed by $p$.
Then we have an isomorphism
\[ Y_S(M/F^{\cyc}) \cong (\Gal(K(F^{\cyc})/F^{\cyc})/p) \ot_{\Zp} M \]
of abelian groups. In particular, $Y_S(M/F^{\cyc})$ is finite if and
only if $\Gal(K(F^{\cyc})/F^{\cyc})$ is finitely generated over
$\Zp$. \el

\bpf Since $pM = 0$, we have $N:= \Hom_{\Zp}(M, \mu_{p^{\infty}}) =
\Hom_{\Zp}(M, \mu_p)$. As the group $G_S(F^{\cyc})$ acts trivially
on $M$ and $\mu_p \sbs F$, it also acts trivially on
$N=\Hom_{\Zp}(M, \mu_p)$. Therefore, we have
\[ R_S(N/F^{\cyc}) = \Hom_{\Zp}(\Gal(K(F^{\cyc})/F^{\cyc})/p, N), \]
noting that every finite prime of $F^{\cyc}$ splits completely in
$K(F^{\cyc})$. On the other hand, one has the following adjunction
isomorphism
\[ \Big(\Gal(K(F^{\cyc})/F^{\cyc})/p \ot_{\Zp}M\Big)^{\vee}
\cong \Hom_{\Zp}(\Gal(K(F^{\cyc})/F^{\cyc})/p, N).\] (Here we are
identifying $M^{\vee}$ with $N$ using the assumption that
$F(\mu_p)=F$ and $pM=0$.) Combining both equalities and taking dual,
we obtain the required isomorphism. The second assertion is
immediate from the first. \epf

\bl \label{fg H2} $($Compare with \cite[Lemma $3.2$]{CS}$)$
$Y_S(T/F^{cyc})$ is finitely generated over $R$ if and only if
$H^2_{S}(F^{cyc}/F, T)$ is finitely generated over $R$. \el

\bpf From the Poitou-Tate sequence, we have the following exact
sequence
\[ 0\lra Y_S(T/F^{\cyc}) \lra H^2_{S}(F^{cyc}/F, T)\lra \Big(\bigoplus_{v\in S}K_v^0(W/F^{\cyc})\Big)^{\vee}. \]
  For each $v\in S$, one has
\[ K_v^0(W/F^{\cyc}) \cong \bigoplus_{w|v} W(F^{cyc}_w).  \]
Since every primes splits finitely in $F^{\cyc}/F$, the direct sum
is a finite sum of cofinitely generated $R$-modules, and so it is a
cofinitely generated $R$-module. (Note that the direct sum needs not
be a finite sum if $p=2$ and $F$ has real primes but our standing
assumption is that if $p=2$, then $F$ has no real primes, so this
situation will not occur under our assumption.) The conclusion of
the lemma then follows. \epf

We can now prove our theorem.

\bpf[Proof of Theorem \ref{main}]
  By Lemma \ref{fg finite}, we may replace $F$, if necessary, so
  that $F$ contains $\mu_{2p}$ (in particular, $F$ has no real primes) and that $G_S(F)$ acts trivially on
  $T/\m T$. By Lemma \ref{fg H2}, it then suffices to show that
  $H^2_S(F^{\cyc}/F, T)$ is finitely generated over $R$. Choose a
  set of generators $x_1,..., x_d$ of $\m$ such that they form a
  basis for $\m/\m^2$. By a repeated application of Lemma \ref{residual
  H2}, we have an isomorphism $H^2_S(F^{\cyc}/F, T)/\m \cong H^2_S(F^{\cyc}/F, T/\m T)$.
  By Nakayama lemma, we are reduced to showing that $H^2_S(F^{\cyc}/F, T/\m T)$ is
  finite. By another application of Lemma \ref{fg H2}, this is equivalent to showing that
  $Y((T/\m T)/ F^{\cyc})$ is finite. This latter assertion is an immediate consequence
  of Lemma \ref{fine selmer and class
  group} and the validity of the Iwasawa $\mu$-conjecture.
\epf

Upon a finer examination of our proof of Theorem \ref{main}, one
actually shows something more precise. Let $\bar{\rho} :
G_{S}(F)\lra \mathrm{Aut}_k(T/\m T)$ be the residual representation
of $\rho$. Denote $F(\mu_{2p}, T/\m T)$ to be $F_S^{\ker\bar{\rho}}
(\mu_{2p})$. Note that this is a finite Galois extension of $F$.

\bt \label{main prop} Let $L$ be a finite extension of $F$ such that
$F(\mu_{2p}, T/\m T)$ is contained in a finite $p$-extension of $L$.
 Then the Iwasawa $\mu$-invariant conjecture holds
for $L$ if and only if $Y_S(T/L^{\cyc})$ is finitely generated over
$R$. \et

\bpf Let $L'$ be a finite $p$-extension of $L$ that contains
$F(\mu_{2p}, T/\m T)$. The proof of Theorem \ref{main} essentially
proved the equivalence over $L'$. To establish the equivalence for
$L$, we need to show that the finite generation property is
preserved in a finite $p$-extension. By \cite[Theorem 3]{Iw}, the
Iwasawa $\mu$-invariant conjecture holds for $L^{\cyc}$ if and only
if the Iwasawa $\mu$-invariant conjecture holds for $L'^{\cyc}$. It
remains to show the same for the case of $Y_S(T/L^{\cyc})$. It is
not difficult, by making use of the commutative diagram in Lemma
\ref{fg finite}, to show that the map
\[ Y_S(T/L'^{\cyc})_{\triangle}\longrightarrow Y_S(T/L^{\cyc})\]
has kernel and cokernel which are finitely generated over $R$. Here
$\triangle=\Gal(L'/L)$. It follows from this observation that
$Y_S(T/L^{\cyc})$ is finitely generated over $R$ if and only if
$Y_S(T/L'^{\cyc})_{\triangle}$ is finitely generated over $R$. Since
$\triangle$ is a $p$-group, $R[\triangle]$ is local with a unique
maximal (two-sided) ideal $\M =\m R[\triangle]+I_{\triangle}$, where
$I_{\triangle}$ is the augmentation ideal (see \cite[Proposition
5.2.16(iii)]{NSW}). It is easy to see from this that
 \[ Y_S(T/L'^{\cyc})/\M \cong Y_S(T/L'^{\cyc})_{\triangle}/\m Y_S(T/L'^{\cyc})_{\triangle}.
  \] Therefore, Nakayama's lemma
for $R$-modules tells us that $Y_S(T/L'_{\cyc})_{\triangle}$ is
finitely generated over $R$ if and only if $Y_S(T/L'^{\cyc})/\M$ is
finite. On the other hand, Nakayama's lemma for
$R[\triangle]$-modules tells us that $Y_S(T/L'^{\cyc})/\M$ is finite
if and only if $Y_S(T/L'^{\cyc})$ is finitely generated over
$R[\triangle]$. But since $\triangle$ is finite, the latter is
equivalent to saying that $Y_S(T/L'^{\cyc})$ is finitely generated
over $R$. Hence we conclude that $Y_S(T/L^{\cyc})$ is finitely
generated over $R$ if and only if $Y_S(T/L'^{\cyc})$ is finitely
generated over $R$. The proof of the theorem is now completed. \epf

One expects that the conjecture is invariant under isogeny. However,
at present, we can only establish this partially in the following
proposition. We will state this proposition in a more general
context.  Suppose that $R'$ is another commutative Noetherian local
ring with maximal ideal $\m'$ and finite residue field of
characteristic $p$, and suppose that $T'$ is a finitely generated
$R'$-module with a continuous $R$-linear $G_S(F)$-action. We then
have the following proposition.

\bp \label{isogeny} Suppose that $F(T/\m T, T'/\m'T', \mu_{2p})$ is
a finite $p$-extension of $F$. Then $Y_S(T/F^{\cyc})$ is finitely
generated over $R$ if and only if $Y_S(T'/F^{\cyc})$ is finitely
generated over $R'$.\ep

\bpf As seen in the proof of Theorem \ref{main prop}, the finite
generation property is preserved in a finite $p$-extension.
Therefore, one is reduced to showing that $Y_S(T/L^{\cyc})$ is
finitely generated over $R$ if and only if $Y_S(T'/L^{\cyc})$ is
finitely generated over $R'$, where $L: = F(T/\m T, T'/\m'T
,\mu_{2p})$. Since $L$ is clearly a finite $p$-extension of $F(T/\m
T,\mu_{2p})$ (resp., $F(T'/\m'T ,\mu_{2p})$), Theorem \ref{main
prop} applies to imply the assertion that the Iwasawa
$\mu$-invariant conjecture holds for $L$ if and only if
$Y_S(T/F^{\cyc})$ is finitely generated over $R$ (resp.,
$Y_S(T'/F^{\cyc})$ is finitely generated over $R'$). The proposition
then follows. \epf

We end the section with the following remarks.

\begin{remarks} (a) In the case when $T$ is the Tate module of an elliptic
curve and the prime $p$ is odd, the field extension $F(\mu_{2p},
T/\m T)$ is precisely $F(E[p])$. One can easily see that Theorem
\ref{main prop} recovers \cite[Theorem 3.4]{CS} and Proposition
\ref{isogeny} recovers the observation made after \cite[Lemma
3.8]{CS}.

\smallskip
 \noindent (b) One can give alternative proofs to
 \cite[Theorem 3]{Jh} and \cite[Theorem 8]{JhS} by appealing to Lemma \ref{fg H2} and Lemma \ref{residual
 H2}.

\smallskip \noindent
(c) The results here can be extended easily to the case when the
ring $R$ is a commutative Noetherian semi-local ring, complete with
respect to its Jacobson radical $J(R)$, and that $R/J(R)$ is a
finite ring of order a power of $p$. In this case, the ring $R$ has
finitely many maximal ideals $\m_1,...,\m_r$ and is isomorphic to
\[ R_{\m_1}\times\cdots\times R_{\m_r}, \]
where each $R_{\m_i}$ is commutative complete Noetherian local with
finite residue field of characteristic $p$ (see \cite[Theorem
8.15]{M}. Of course, there is still some work to be done after
applying the said theorem, namely, one still needs to show that each
$R_{\m_i}$ is $\m_i$-adic complete but this can be easily verified).
Now every $R[G_S(F)]$-module $T$ decomposes canonically as\[
T_{\m_1}\times\cdots\times T_{\m_r},
\] and it is easy to see that the decomposition is compatible with
the Galois action. The results in this section can then be applied
to each $R_{\m_i}$-module $T_{\m_i}$.

\end{remarks}

\section{Ranks of Iwasawa modules} \label{Ranks of Iwasawa
modules}

In this section, we establish some algebraic preliminaries to
facilitate further discussion of the fine Selmer groups. We will
prove certain formulas on the rank of modules over a completed group
ring. Much of the materials considered here originates from
\cite{BH, HO, Ho}. As a start, we shall prove a general lemma. Let
$\La$ be a (not neccessarily commutative) Noetherian ring which has
no zero divisors. Then it admits a skew field of fractions $K(\La)$
which is flat over $\La$ (see \cite[Chapters 6 and 10]{GW} or
\cite[Chapter 4, \S 9 and \S 10]{Lam}). If $M$ is a finitely
generated $\La$-module, we define the $\La$-rank of $M$ to be
$$ \rank_{\La}M  = \dim_{K(\La)} K(\La)\ot_{\La}M. $$
Clearly, one has $\rank_{\La}M =0 $ if and only if
$K(\La)\ot_{\La}M=0$.
The following lemma will be useful in the discussion in this
section.
 \bl \label{La rank}  Let $\La$ be a Noetherian ring which has no
zero divisors. Suppose $\Omega$ is a quotient of $\La$ such that it
is also a Noetherian ring which has no zero divisors. Let $M$ be a
finitely generated $\La$-module which has a finite free resolution
of finite length. Then we have
 $$\rank_{\La}M = \sum_{i\geq 0}(-1)^i\rank_{\Omega}\Tor^{\La}_i(\Omega, M).
 $$\el

 \bpf (Compare with proof of \cite[Theorem 1.1]{Ho}) Let $$ 0\lra \La^{n_d}\lra \cdots \lra \La^{n_0}\lra M $$ be a
resolution of $M$ which exists by the assumptions of the lemma. Then
the groups $\Tor^{\La}_i(\Omega, M)$ can be computed by the homology
of the complex
$$  \Omega^{n_d}\lra \cdots \lra \Omega^{n_0}.$$
This in turn implies that each $\Tor^{\La}_i(\Omega, M)$ is finitely
generated over $\Omega$ and the
 sum on the right hand side is a finite sum whose value coincides
 with
 $$ \sum_{i= 0}^d (-1)^i n_i.$$
But this latter quantity is precisely $\rank_{\La}M$.
 \epf

We record another lemma. Write $M^+ = \Hom_{\La}(M,\La)$.

\bl \label{torsion is Hom zero} Let $\La$ be a Auslander regular
ring (see \cite[Definition 3.3]{V1}) with no zero divisors. Let $M$
be a finitely generated $\La$-module. Then the following are
equivalent.

\smallskip $(a)$ The canonical map
 $\phi: M\lra M^{++}$ is zero.

\smallskip $(b)$ $K(\La)\ot_{\La}M =0$, where $K(\La)$ is the skew
field of $\La$.

\smallskip $(c)$ $\Hom_{\La}(M,\La)= 0$.
 \el

\bpf The equivalence of (a) and (c) follows from \cite[Remark
3.7]{V1}. Suppose that $K(\La)\ot_{\La}M = 0$. Let $f\in
\Hom_{\La}(M,\La)$ and $x\in M$. Then since $K(\La)\ot_{\La}M = 0$,
there exists $\la\in\La\setminus\{0\}$ such that $\la x=0$. This in
turn implies that $\la f(x) = f(\la x) = 0$. Since $\La$ has no zero
divisor, we have $f(x)=0$. This shows that $\Hom_{\La}(M,\La)= 0$
and the implication (b)$\Rightarrow$(c). Conversely, suppose that
$\Hom_{\La}(M,\La)= 0$. By \cite[Proposition 2.5]{V1} and the
Auslander condition, the canonical map
 $\phi: M\lra M^{++}$
has kernel and cokernel which are $R\ps{H}$-torsion. Therefore,
$\phi$ induces an isomorphism
$$K(\La)\ot_{\La}M \stackrel{\sim}{\lra} K(\La)\ot_{\La}M^{++}.$$
Now if $\phi =0$, then it will follow immediately that
$K(\La)\ot_{\La}M = 0$. This establishes (a)$\Rightarrow$(b). \epf

We now apply the above discussion to the context of a completed
group ring. Let $R$ be a complete regular local ring with finite
residue field of characteristic $p$, where $p$ is a prime. Let $H$
be a compact pro-$p$ $p$-adic Lie group without $p$-torsion. It is
well known that $R\llbracket H\rrbracket$ is a Auslander regular
ring (see \cite{V1}, see also Theorem \ref{compact Auslander
regular}). In particular, the ring $R\ps{H}$ is Noetherian local and
has finite projective dimension. Therefore, it follows that every
finitely generated $R\ps{H}$-module admits a finite free resolution
of finite length. In the case that either the ring $R$ has
characteristic zero or $H$ is a uniform pro-$p$ group, the ring
$R\ps{H}$ has no zero divisors (see Theorem \ref{compact Auslander
regular} and remarks after it), and therefore admits a skew field
which enable one to define the notion of a rank as above. Now
suppose that $R$ has characteristic $p$ and $H$ is a compact pro-$p$
$p$-adic Lie group without $p$-torsion. We then define the
$R\ps{H}$-rank of a finitely generated $R\ps{H}$-module $M$ by
 $$ \rank_{R\ps{H}} M =
 \displaystyle\frac{\rank_{R\ps{H_0}}M}{|H:H_0|}, $$
where $H_0$ is an open normal uniform pro-$p$ subgroup of $H$. We
will see below that this is integral and independent of the choice
of $H_0$.

\bl \label{Tor rank} Let $H$ be a compact pro-$p$ $p$-adic Lie group
without $p$-torsion. Let $M$ be a finitely generated
$R\ps{H}$-module. Then $H_i(H,M)$ is finitely generated over $R$ for
each $i$ and we have the equality
 $$  \rank_{R\ps{H}}M =  \displaystyle\sum_{i\geq 0} (-1)^i\rank_{R}H_i(H,M)\\
     $$ In particular, the definition of $R\ps{H}$-rank is well-defined and integral.
     \el

\bpf Suppose first that either the ring $R$ has characteristic zero
or $H$ is a uniform pro-$p$ group. Then the conclusion follows from
applying Lemma \ref{La rank} (taking $\La =R\ps{H}$ and $\Omega =
R$) and observing that $H_i(H,M)\cong \Tor^{R\ps{H}}_i(R,M)$ for all
$i$. Now suppose that $R$ has characteristic $p$ and $H$ is a
compact pro-$p$ $p$-adic Lie group without $p$-torsion. Fix a finite
free $R\ps{H}$-resolution
$$ 0\lra R\ps{H}^{n_d} \lra \cdots \lra R\ps{H}^{n_0}\lra M$$ of $M$. Then the groups $H_i(H,M) =
\Tor^{R\ps{H}}_i(R, M)$ can be computed by the homology of the
complex
$$  R^{n_d}\lra \cdots \lra R^{n_0}.$$
This in turn implies that
 $$ \sum_{i=0}^d (-1)^i\rank_{R}H_i(H,M)= \sum_{i= 0}^d (-1)^i n_i.$$ Let $H_0$ be any open normal uniform
pro-$p$ subgroup of $H$. The above free $R\ps{H}$-resolution is also
a free $R\ps{H_0}$-resolution for $M$. Therefore, the groups
$H_i(H_0,M) = \Tor^{R\ps{H_0}}_i(R, M)$ can be computed by the
homology of the complex
$$  R^{|H:H_0|n_d}\lra \cdots \lra R^{|H:H_0|n_0}$$ which gives
 $$ \sum_{i=0}^d (-1)^i\rank_{R}H_i(H_0,M)= |H:H_0|\sum_{i= 0}^d (-1)^i n_i.$$
On the other hand, the sum on the left is precisely
$\rank_{R\ps{H_0}}M$ by an application of Lemma \ref{La rank}
(taking $\La =R\ps{H_0}$ and $\Omega = R$). Hence, we have
$$ \frac{\rank_{R\ps{H_0}}M}{|H:H_0|} = \sum_{i= 0}^d (-1)^i n_i.$$
The sum on the right is clearly integral and independent of $H_0$.
Therefore, we have proved the proposition.
 \epf

We may now define the notion of torsion modules over $R\ps{H}$ via
the following lemma.

\bl \label{torsion is Hom zero RG} Let $R$ be a regular local ring
with finite residue field of characteristic $p$. Let $H$ be a
compact pro-$p$ $p$-adic Lie group without $p$-torsion. Let $M$ be a
finitely generated $R\ps{H}$-module. Then $\rank_{R\ps{H}}M =0$ if
and only if $\Hom_{R\ps{H}}(M,R\ps{H})=0$.\el

\bpf If the regular local ring $R$ has characteristic zero or $H$ is
a uniform pro-$p$ group, then this lemma follows from lemma
\ref{torsion is Hom zero}. For the exceptional case, we let $H_0$ be
an open uniform normal subgroup of $H$.  By the $R$-analog of
\cite[Proposition 5.4.17]{NSW}, we have
$$ \Hom_{R\ps{H}}(M, R\ps{H}) \cong \Hom_{R\ps{H_0}}(N,
R\ps{H_0}). $$ Therefore, it follows that
$\Hom_{R\ps{H}}(M,R\ps{H})=0$ if and only if
$\Hom_{R\ps{H_0}}(M,R\ps{H_0})=0$. On the other hand, it is clear
from the definition that $\rank_{R\ps{H}}M =0$ if and only if
$\rank_{R\ps{H_0}}M =0$. Thus, we may then apply the above
discussion to obtain the required equivalence. \epf

In view of the above lemma, we will say that a finitely generated
$R\ps{H}$-module $M$ is \textit{torsion} if either of the two
equivalent statements hold. The following lemma is a relative
version of Lemma \ref{Tor rank}.

\bl \label{relative rank} Let $H$ be a compact pro-$p$ $p$-adic Lie
group without $p$-torsion. Let $U$ be a closed normal subgroup of
$H$ such that $H/U$ is also a compact pro-$p$ $p$-adic Lie group
without $p$-torsion. Let $M$ be a finitely generated
$R\ps{H}$-module. Then $H_i(U,M)$ is finitely generated over
$R\ps{H/U}$ for each $i$ and we have the equality
 $$ \ba{rl} \rank_{R\ps{H}}M = & \displaystyle\sum_{i\geq 0} (-1)^i\rank_{R\ps{H/U}}H_i(U,M)\\
    =& \rank_{R\ps{H/U}}M_U + \displaystyle\sum_{i\geq 1}
  (-1)^i\rank_{R\ps{H/U}}H_i(U,M). \ea
 $$ \el

\bpf If the ring $R$ has characteristic zero or $H$ is a uniform
pro-$p$ group, the conclusion follows from applying Lemma \ref{La
rank} (taking $\La =R\ps{H}$ and $\Omega = R\ps{H/U}$) and observing
that $H_i(U,M)\cong \Tor^{R\ps{H}}_i(R\ps{H/U},M)$ for all $i$. In
general, one has the equality
 $$ \ba{rl} \rank_{R\ps{H}}M = & \displaystyle\sum_{i\geq 0} (-1)^i\rank_{R}H_i(H,M)\\
    =&  \displaystyle\sum_{i,j \geq 0}
  (-1)^{i+j} \rank_{R}H_i(H/U, H_j(U,M)) \\
  =&  \displaystyle\sum_{j \geq 0}
  (-1)^{j} \rank_{R\ps{H/U}}H_i(U,M), \ea
 $$
where the first and third equality follows from Lemma \ref{Tor rank}
and the second equality is a consequence of the spectral sequence
 $$ H_i(H/U, H_j(U, M))\Longrightarrow H_{i+j}(H,M). $$
\epf

It is naively immediate from the equality in Lemma \ref{relative
rank} that if
$$\rank_{R\ps{H}}M = \displaystyle\sum_{i\geq 1}
(-1)^i\rank_{R\ps{H/U}}H_i(U,M) =0, $$ then $\rank_{R\ps{H/U}}M_U
=0$. The converse is true if we assume further that $U$ is solvable,
and this is the content of the next theorem.

\bt \label{torsion main theorem/lemma}
 Let $H$ be a compact pro-$p$ $p$-adic Lie
group without $p$-torsion. Let $U$ be a closed normal subgroup of
$H$ such that $U$ is solvable and $H/U$ is a compact pro-$p$
$p$-adic Lie group without $p$-torsion. Let $M$ be a finitely
generated $R\ps{H}$-module. Then $M_U$ is a torsion
$R\ps{H/U}$-module if and only if $M$ is a torsion $R\ps{H}$-module
and
 $\displaystyle\sum_{i\geq 1}
(-1)^i\rank_{R\ps{H/U}}H_i(U,M) =0$.\et

\bpf The above discussion already establishes one direction. Now
suppose that $M_U$ is a torsion $R\ps{H/U}$-module. The required
conclusion will follow once we show that $M$ is a torsion
$R\ps{H}$-module. This will follow from the next lemma which is a
slight refinement of the final theorem in \cite{BH} and \cite[Lemma
2.6]{HO}.
  \epf

\bl \label{torsion lemma} Let $H$ be a compact $p$-adic Lie group
without $p$-torsion. Suppose that $N$ is a closed normal subgroup of
$H$ with the property that $N$ is a solvable uniform pro-$p$ group
and $H/N$ has no $p$-torsions. Let $M$ be a finitely generated
$R\ps{H}$-module. If $M_{N}$ is a torsion $R\ps{H/N}$-module, then
$M$ is a torsion $R\ps{H}$-module.  \el

To prove this lemma, we need two more lemmas.

\bl \label{uniform lemma} Let $H$ be a uniform pro-$p$ group, and
let $N$ be a closed normal subgroup of $H$ with the property that
$N$ is a solvable uniform pro-$p$ group and $H/N$ has no
$p$-torsion. Then there exists a closed normal subgroup $N_0$ of $H$
satisfying all the following properties.

\smallskip
$(i)$  $N_0 \cong \Zp^r$ for some $r>0$.

\smallskip
$(ii)$ $H/N_0$ is uniform with $\dim H/N_0<\dim H$.

\smallskip
$(iii)$ $N_0\sbs N$.\el

\bpf If $N$ is abelian, then one may take $N_0 = N$. Now if $N$ is
not abelian, we write $N^{0}= N$ and $N^{(n+1)} =
\overline{[N^{(n)}, N^{(n)}]}$. Then $N^{(m+1)} = 0$ but
$N^{(m)}\neq 0$ for some $m\geq 1$. Note that $N^{(m)}$ is abelian.
Set
\[ N_0 : = \{h\in H ~|~ h^{p^j}\in N^{(m)}~\mbox{for~some}~j\} \]
The proof of statement (3) of the first proposition in \cite[\S
4]{BH} shows that $N_0$ satisfies (i) and (ii). To see that $N_0$
satisfies (iii), one applies a similar argument as in the last
paragraph of the proof of \cite[Lemma 2.6]{HO}.\epf

\bl \label{torsion lemma2} Let $H$ be a uniform pro-$p$ group, and
let $N$ be a closed normal subgroup of $H$ with the property that
$N\cong \Zp^r$ and $H/N$ has no $p$-torsion.  Let $M$ be a finitely
generated $R\ps{H}$-module. If $M_{N}$ is a torsion
$R\ps{H/N}$-module, then $M$ is a torsion $R\ps{H}$-module. \el

\bpf We prove this lemma by the method of contradiction, following
the argument given in the last theorem in \cite{BH}. Suppose that
$M$ is a finitely generated $R\ps{H}$-module with $R\ps{H}$-rank
$s>0$. Recall that $R\ps{H}$ is Auslander regular and has no
zero-divisors (cf. Theorem \ref{compact Auslander regular}). We
shall first show that there is map $M \lra R\ps{H}^s$ with
$R\ps{H}$-torsion kernel and cokernel. Let $K(H)$ denote the skew
field of $R\ps{H}$. Write $M^+ = \Hom_{R\ps{H}}(M, R\ps{H})$. Then
by \cite[Proposition 2.5]{V1}, there is a map $M\lra M^{++}$ with
$R\ps{H}$-torsion kernel and cokernel (the torsionness comes from
the definition of Auslander regularity). Choose $f_1,..., f_s\in
M^+$ such that they form a basis for $K(H)\ot_{R\ps{H}}M^+$. Then
these elements give rise to a map $R\ps{H}^s\lra M^+$ which clearly
has $R\ps{H}$-torsion kernel and cokernel. Taking $R\ps{H}$-dual, we
obtain a map $M^{++}\lra R\ps{H}^s$ with $R\ps{H}$-torsion kernel
and cokernel. Combining this with the above canonical map, we obtain
the required map.

Projecting onto any (fixed) factor of $R\ps{H}^s$, we obtain a
$R\ps{H}$-homomorphism $\phi:M\lra R\ps{H}$. It is not difficult to
see that $\phi$ is nontrivial. Denote $I(N)$ to be augmentation
ideal of the ring $R\ps{N}$, and denote $J$ to be the two-sided
ideal $I(N)R\ps{H} = R\ps{H}I(N)$. Since $J$ is a closed ideal of
$R\ps{H}$, we have $\cap_{n\geq 0} J^n = 0$. Therefore, we can find
an $m$ such that $\phi(M)\sbs J^m$ and $\phi(M)\nsubseteq J^{m+1}$.
Then $M' = (\phi(M)+J^{m+1})/J^{m+1}$ is a nontrivial submodule of
$J^m/J^{m+1}$. We also note that $N$ acts trivially on $M'$ and
$J^m/J^{m+1}$. Hence we have $M' = (M')_N$ being a quotient of
$M_N$. On the other hand, by a similar argument to that in the last
paragraph of \cite{BH}, we have that $J^m/J^{m+1}$ is a free
$R\ps{H/N}$-module with positive $R\ps{H/N}$-rank. Since $M'$ is a
nontrivial submodule of a free $R\ps{H/N}$-module, $M'$ must also
have positive $R\ps{H/N}$-rank which in turn implies that $M_N$ has
positive $R\ps{H/N}$-rank and this contradicts the hypothesis of the
lemma.
 \epf

We can now prove Lemma \ref{torsion lemma}.

\bpf[Proof of Lemma \ref{torsion lemma}]
 We proceed by induction on the dimension of $H$. When $\dim H=1$,
 the assertion of the lemma can be deduced from the classical result
 of Iwasawa. Now suppose $\dim H>1$. If $N=\Zp^r$, then we are done
 by Lemma \ref{torsion lemma2}. Else by Lemma \ref{uniform lemma},
 we can find a closed normal subgroup $N_0$ of $H$ such that (i) $N_0 \cong \Zp^r$ for some
 $r>0$, (ii) $H/N_0$ is uniform with $\dim H/N_0<\dim H$ and (iii)
$N_0\sbs N$. Viewing $M_N = (M_{N_0})_{N/N_0}$, we may apply our
induction hypothesis to deduce that $M_{N_0}$ is a torsion
$R\ps{H/N_0}$-module. Now applying Lemma \ref{torsion lemma2} again,
we obtain the required conclusion. \epf

The following is an immediate corollary of the above. Alternatively,
one can prove this directly by the rank calculation.

\bc Let $H$ be a compact pro-$p$ $p$-adic Lie group without
$p$-torsion. Let $U$ be a closed normal subgroup of $H$ such that
$U\cong \Zp$ and $H/U$ is a compact pro-$p$ $p$-adic Lie group
without $p$-torsion. Let $M$ be a finitely generated
$R\ps{H}$-module. Then $M_U$ is a torsion $R\ps{H/U}$-module if and
only if $M$ is a torsion $R\ps{H}$-module and $H_1(U,M)$ is a
torsion $R\ps{H/U}$-module.\ec

We mention another corollary of Lemma \ref{torsion lemma} which can
be proved via a similar argument as in \cite[Lemma 2.5]{HO}.

\bc \label{torsion lemma corollary} Let $H$ be a compact $p$-adic
Lie group without $p$-torsion. Suppose that $N$ is a closed normal
subgroup of $H$ with the property that $N$ is a solvable uniform
pro-$p$ group and $H/N$ has no $p$-torsions. Let $M$ be a finitely
generated $R\ps{H}$-module. Then we have $$\rank_{R\ps{G}}M \leq
\rank_{R\ps{H/N}}M_N.$$ \ec

We now mention another variant of Howson's result which can be
viewed as a slight generalization of \cite[Corollary 1.10]{Ho} (see
also \cite[Lemma 5]{Jh}).

\bp \label{torsion x}  Let $H$ be a compact $p$-adic Lie group
without $p$-torsion, and let $x\in R$ be a nonzero element such that
$\bar{R}: = R/xR$ is also a complete regular local ring. Let $M$ be
a finitely generated $R\ps{H}$-module. Then
  $$ \rank_{\bar{R}\ps{H}} M/xM = \rank_{\bar{R}\ps{H}}M[x] +
  \rank_{R\ps{H}}M,
 $$ where $M[x]$ is the submodule of $M$ killed by $x$. \ep

\bpf Using the free resolution of $\bar{R}\llbracket G\rrbracket$
\[ 0\lra R\ps{G} \stackrel{x}{\lra} R\ps{G} \lra \bar{R}\ps{G}\lra 0,\]
 one can calculate that \[  \Tor^i_{R\ps{H}} (\bar{R}\ps{H}, M) =
\begin{cases} M/xM & \text{\mbox{if} $i=0$},\\
  M[x] & \text{\mbox{if} $i=1$},\\
  0 & \text{\mbox{otherwise}.}
\end{cases} \] The conclusion is now immediate from the above calculations and Lemma \ref{La rank}. \epf

\bc \label{torsion x corollary} Retaining the above assumptions. Let
$M$ be a finitely generated $R\ps{H}$-module. Then
  $M/xM$ is a torsion $\bar{R}\ps{H}$-module if and only if $M[x]$ is a torsion $\bar{R}\ps{H}$-module
  and $M$ is a torsion $R\ps{H}$-module. \ec

\section{On pseudo-nullity of fine Selmer groups} \label{pseudo-nullity section}

In this section, we will investigate certain properties of the fine
Selmer groups over $p$-adic Lie extensions of dimension strictly
larger than one. As in the previous section, we shall assume that
our ring $R$ is a complete regular local ring with finite residue
field of characteristic $p$. Let $G$ be a compact pro-$p$ $p$-adic
Lie group without $p$-torsion. Then $R\ps{G}$ is a Auslander regular
ring. We say that a finitely generated torsion $R\ps{G}$-module is
\textit{pseudo-null} if $ \Ext_{R\ps{G}}^{1}(M, R\ps{G}) = 0$. The
following fundamental lemma will be crucial in our discussion.

\bl \label{pseudo-null torsion} Let $R$ be a complete regular local
ring with finite residue field of characteristic $p$ and let $G$ be
a compact pro-$p$ $p$-adic Lie group without $p$-torsion. Suppose
that $H$ is a closed normal subgroup of $G$ with $G/H \cong \Zp$.
Let $M$ be a compact $R\ps{G}$-module which is finitely generated
over $R\ps{H}$. Then $M$ is a pseudo-null $R\ps{G}$-module if and
only if $M$ is a torsion $R\ps{H}$-module. \el

\bpf This essentially follows from \cite{V2}, which we shall explain
(see also \cite[Lemma 3.1]{HS}). By the standard theory of compact
$p$-adic Lie groups (for instances, see \cite{DSMS}), one can find
open subgroups $H_0$ of $H$ and $G_0$ of $G$ such that $H_0$ and
$G_0$ are uniform pro-$p$ groups and $G_0/H_0 \cong \Zp$. By
\cite[Proposition 5.4.17]{NSW} (the same conclusion with a similar
proof holds if we replace $\Zp$ by $R$), we have
\[ \Ext_{R\ps{G}}^{i}(N, R\ps{G}) \cong
\Ext_{R\ps{G_0}}^{i}(N, R\ps{G_0})\] for any $R\ps{G}$-module $N$
and all $i$, and a similar statement holds for $H$ and $H_0$.
Therefore, we are reduced to showing the lemma under the assumptions
that $G$ and $H$ are uniform pro-$p$ groups.

We now write $\Ga = G/H$. There is a natural group homomorphism
$\phi: \Ga \lra \mathrm{Aut}(H)$. Suppose that
\[ \im \phi \sbs \{ f\in \mathrm{Aut}(H)~| ~f(h)h^{-1} \in H^p~\mbox{for all}~h\in H\}. \hspace{0.4in} (\ast)\]
Then the conclusion of the lemma follows from an $R$-analog of
\cite[Example 2.3]{V2} and \cite[Proposition 5.4]{V2} in this
instance. In general, since $H$ is a uniform pro-$p$ group by our
assumption, $H^p$ is an open characteristic subgroup of $H$.
Therefore, the map $\phi$ induces a continuous group homomorphism
$\bar{\phi}: \Ga \lra \mathrm{Aut}(H/H^p)$. Since $H/H^p$ is finite,
we have $\Ga_1:=\ker\bar{\phi} \cong \Zp$. Let $G_1$ be the open
normal subgroup of $G$ containing $H$ such that $G_1/H = \Ga_1$. By
another application of the $R$-analog of \cite[Proposition
5.4.17]{NSW}, we are reduced to showing that $M$ is a pseudo-null
$R\llbracket G_1\rrbracket$-module if and only if $M$ is a torsion
$R\llbracket H\rrbracket$-module. However, in this case, the natural
map $\phi_1 : \Ga_1 \lra \mathrm{Aut}(H)$ clearly satisfies $(\ast)$
by our choice of $\Ga_1$, and so \cite[Example 2.3, Proposition
5.4]{V2} can be applied. \epf

\begin{remark} In the case when $H\cong \Zp$, the lemma follows directly from
the $R$-analog of \cite[Example 2.3]{V2} and \cite[Proposition
5.4]{V2}. \end{remark}

We return to arithmetic. As before, $p$ will denote a prime. For the
remainder of the paper, we shall assume further that the number
field $F$ has no real primes when $p=2$.
 Following \cite{CS}, we say that $F_{\infty}$ is a
$S$-admissible $p$-adic extension of $F$ if (i) $\Gal(F_{\infty}/F)$
is compact pro-$p$ $p$-adic Lie group without $p$-torsion, (ii)
$F_{\infty}$ contains $F^{\cyc}$ and (iii) $F_{\infty}$ is contained
in $F_S$. Write $G=\Gal(F_{\infty}/F)$ and
$H=\Gal(F_{\infty}/F^{\cyc})$. From now on, we shall assume that our
Galois module $T$ is a free $R$-module of finite rank with a
continuous $R$-linear $G_S(F)$-action. The following lemma is also
considered in \cite[Lemma 3.2]{CS}.

\bl \label{fg La H} Let $F_{\infty}$ be a $S$-admissible $p$-adic
Lie extension of $F$. Then the following statements are equivalent.

\smallskip
$(a)$ $Y_S(T/F^{\mathrm{cyc}})$ is a finitely generated $R$-module.

\smallskip
$(b)$ $H^2_S(F^{\cyc}/F, T)$ is a finitely generated $R$-module.

\smallskip
$(c)$ $Y_S(T/F_{\infty})$ is a finitely generated $R\ps{H}$-module.

\smallskip
$(d)$  $H^2_S(F_{\infty}/F, T)$ is a finitely generated
$R\ps{H}$-module.\el

\bpf The equivalence of (a) and (b) is shown in Lemma \ref{fg H2}.
The equivalence of (c) and (d) can be shown by a similar argument.
The equivalence of (b) and (d) follows from Lemma \ref{descent H2}
and Nakayama lemma. \epf

The following question has been studied by many before.

\medskip
\noindent \textbf{Question B:}  Let $F_{\infty}$ be a $S$-admissible
$p$-adic Lie extension of $F$ of dimension $>1$, and suppose that
$Y_S(T/F_{\infty})$ is a finitely generated $R\ps{H}$-module. Is
$Y_S(T/F_{\infty})$ a pseudo-null $R\ps{G}$-module, or equivalently
a torsion $R\ps{H}$-module?

\medskip
This is precisely \cite[Conjecture B]{CS} when $T$ is the Tate
module of an elliptic curve. In this context, the conjecture has
also been studied and verified for some elliptic curves in \cite{Bh}
and \cite{Oc}. When $T$ is the $R(1)$-dual of the Galois
representation attached to a normalized eigenform ordinary at $p$,
this is \cite[Conjecture B]{Jh}. In the case when $T$ is the
$R(1)$-dual of the Galois representation coming from a $\La$-adic
form, this is \cite[Conjecture 2]{Jh}. We note that in the case when
$T$ is the Tate module of all the $p$-power roots of unity, the dual
fine Selmer group is precisely $\Gal(K(F_{\infty})/F_{\infty})$,
where $K(F_{\infty})$ is the maximal unramified pro-$p$ extension of
$F_{\infty}$ at which every prime of $F_{\infty}$ above $p$ splits
completely, as defined in Section \ref{Cyclotomic Zp extension}. In
this case, Hachimori and Sharifi \cite{HS} has constructed a class
of admissible $p$-adic Lie extension $F_{\infty}$ of $F$ of
dimension $>1$ such that $\Gal(K(F_{\infty})/F_{\infty})$ is not
pseudo-null. Despite so, they have speculated that the
pseudo-nullity condition should hold for admissible $p$-adic
extensions ``coming from algebraic geometry'' (see \cite[Question
1.3]{HS} for details, and see also \cite[Conjecture 7.6]{Sh1} for a
related
assertion and \cite{Sh2} for positive results in this direction). 


Before continuing our discussion, we introduce the following
hypothesis on our admissible extension $F_{\infty}$.

\medskip
\noindent \textbf{(Dim$_{S}$):}  For each $v\in S$, the
decomposition group of $G$ at $v$, denoted by $G_v$, has dimension
$\geq 2$.

\bl \label{pesudo-null H2} Let $F_{\infty}$ be a $S$-admissible
$p$-adic Lie extension of $F$ and assume that $Y_S(T/F_{\infty})$ is
a finitely generated $R\llbracket H\rrbracket$-module. Then the
following statements hold.

\smallskip
$(a)$ If $H^2_{S}(F_{\infty}/F, T)$ is a pseudo-null $R\llbracket
G\rrbracket$-module, so is $Y_S(T/F_{\infty})$.

\smallskip
$(b)$ Suppose that $F_{\infty}$ satisfies $(\mathbf{Dim}_{S})$. Then
if $Y_S(T/F_{\infty})$ is a pseudo-null $R\llbracket
G\rrbracket$-module, so is $H^2_{S}(F_{\infty}/F, T)$. \el

\bpf
  From the Poitou-Tate sequence, we have the following exact
sequence
\[ 0\lra Y_S(T/F_{\infty}) \lra H^2_{S}(F_{\infty}/F, T)\lra \Big(\bigoplus_{v\in S}K_v^0(W/F_{\infty})\Big)^{\vee}. \]
 Statement (a) is then immediate. It remains to show that statement (b) holds.
  For each $v\in S$, fix a prime $w$ of $F_{\infty}$ above $v$. By abuse of notation, we denote the prime of $F^{\cyc}$
  below $w$ by $w$. Write $H_w$ to be the decomposition group of $H$ at $w$. Then one
  sees that $K_v^0(W/F_{\infty})^{\vee}$ is isomorphic to a finite sum
  of terms of the form
\[  R\llbracket H\rrbracket\ot_{R\llbracket H_w \rrbracket}W(F_{\infty,w})^{\vee}
\] which is clearly finitely generated over $R\ps{H}$.
 The assertion of statement (b) will follow once we show that the above term is a torsion $R\llbracket
 H\rrbracket$-module.
 By the assumption $(\mathbf{Dim}_{S})$, the group $H_w$ has dimension
$\geq 1$. It is then easy to see that $W(F_{\infty,w})^{\vee}$ is a
finitely generated torsion $R\llbracket
H_w\rrbracket$-module
. The required conclusion then follows from observing that
\[ \Hom_{R\llbracket
H\rrbracket}\big(R\llbracket H\rrbracket\ot_{R\llbracket H_w
\rrbracket}W(F_{\infty,w})^{\vee}, R\llbracket H \rrbracket\big)=
R\llbracket H\rrbracket\ot_{R\llbracket H_w
\rrbracket}\Hom_{R\llbracket
H_w\rrbracket}\big(W(F_{\infty,w})^{\vee}, R\llbracket H_w
\rrbracket\big) =0. \]
 \epf

We can now state the first main result of this section which is a
slight refinement of the implication $(3)\Rightarrow (1)$ in
\cite[Theorem 10]{Jh}.

\bt \label{pseudo-null main}
 Let $F_{\infty}$ be a $S$-admissible
$p$-adic Lie extension of $F$ and assume that $Y_S(T/F_{\infty})$ is
a finitely generated $R\llbracket H\rrbracket$-module.  Suppose that
there exists a prime ideal $\mathfrak{p}$ of $R$ such that the ring
$R/\mathfrak{p}$ is also regular local. Suppose also that
$F_{\infty}$ satisfies $(\mathbf{Dim}_{S})$.  If
$Y_S\big((T/\mathfrak{p}T)/F_{\infty}\big)$ is a pseudo-null
$R/\mathfrak{p}\llbracket G\rrbracket$-module, then
$Y_S(T/F_{\infty})$ is a pseudo-null $R\llbracket
G\rrbracket$-module. \et

\medskip
Before proving the theorem, we first prove the following lemma.

\bl \label{pseudo-null residual2}
  Let $F_{\infty}$ be a $S$-admissible
$p$-adic Lie extension of $F$. Suppose that $Y_S(T/F_{\infty})$ is a
finitely generated $R\llbracket H\rrbracket$-module.  Let $x$ be a
nonzero and nonunital element of $R$, and suppose that $\bar{R} :=
R/xR$ is also regular local. Write $\bar{T} = T/xT$. Then
$H^2_{S}(F_{\infty}/F, \bar{T})$ is a pseudo-null $\bar{R}\llbracket
G\rrbracket$-module if and only if $H^2_{S}(F_{\infty}/F, T)$ is a
pseudo-null $R\llbracket G\rrbracket$-module and
$H^2_{S}(F_{\infty}/F, T)[x]$ is a pseudo-null $\bar{R}\llbracket
G\rrbracket$-module. \el

\bpf By Lemma \ref{residual H2}, we have an
 isomorphism
 \[ H^2_{S}(F_{\infty}/F, T)/x \cong H^2_{S}(F_{\infty}/F, \bar{T})\]
of $\bar{R}\llbracket G\rrbracket$-modules.
 The lemma is now immediate from
Lemma \ref{pseudo-null torsion} and Corollary \ref{torsion x
corollary}.
 \epf

\medskip
We now give the proof of Theorem \ref{pseudo-null main}.

\bpf[Proof of Theorem \ref{pseudo-null main}] By Lemma
\ref{pesudo-null H2}(a), it suffices to show that
$H^2_{S}(F_{\infty}/F, T)$ is a pseudo-null $R\llbracket
G\rrbracket$-module. By the theory of regular local rings (for
instances, see \cite[\S 14]{M}), one can find a set of generators
$x_1,..., x_r$ of $\mathfrak{p}$ such that each intermediate ring
$R/(x_1,...,x_i)$ is also regular local for $i=1,.., r$. By a
repeated application of Lemma \ref{pseudo-null residual2}, we are
reduced to showing that $H^2_S(F_{\infty}/F,T/\mathfrak{p}T)$ is a
pseudo-null $R/\mathfrak{p}\llbracket G\rrbracket$-module which will
follow from the hypothesis of the theorem and Lemma \ref{pesudo-null
H2}(b). \epf

We record an immediate corollary of Theorem \ref{pseudo-null main}.

\bc \label{pseudo-null corollary}
 Let $F_{\infty}$ be a $S$-admissible
$p$-adic Lie extension of $F$ and assume that $Y_S(T/F_{\infty})$ is
a finitely generated $R\llbracket H\rrbracket$-module. Suppose that
$(\mathbf{Dim}_{S})$ is satisfied.  If
$Y_S\big((T/\mathfrak{m}T)/F_{\infty}\big)$ is a pseudo-null
$k\llbracket G\rrbracket$-module, then $Y_S(T/F_{\infty})$ is a
pseudo-null $R\llbracket G\rrbracket$-module. \ec

We now prove the following descent result for pseudo-nullity.

\bt \label{pseudo-null descent} Let $F_{\infty}$ be a $S$-admissible
$p$-adic Lie extension of $F$ and assume that $Y_S(T/F^{\cyc})$ is a
finitely generated $R$-module. Suppose that $F_{\infty}'$ is another
$S$-admissible $p$-adic Lie extension of $F$ which satisfies the
following properties.

\smallskip
$(i)$  $F'_{\infty}$ is contained in $F_{\infty}$.

\smallskip
$(ii)$ $F_{\infty}'$ satisfies $(\mathbf{Dim}_{S})$.

\smallskip
$(iii)$ The group $N:=\Gal(F_{\infty}/F_{\infty}')$ is a solvable
uniform pro-$p$ group and $H/N$ has no $p$-torsion.

\smallskip
Then $Y_S(T/F_{\infty}')$ is a pseudo-null $R\llbracket
\Gal(F_{\infty}'/F)\rrbracket$-module if and only if
$Y_S(T/F_{\infty})$ is a pseudo-null $R\llbracket
G\rrbracket$-module and $\displaystyle\sum_{i\geq 1}
(-1)^i\rank_{R\ps{H/N}}H_i\big(N, H^2_S(F_{\infty}/F, T)\big) =0$.
\et

\begin{remark} Note that we do not require $F_{\infty}'$ and
$F_{\infty}$ to be solvable extensions of $F$.
\end{remark}

\bpf By Lemmas \ref{descent H2}, \ref{pseudo-null torsion} and
\ref{pesudo-null H2}, it is equivalent to showing that
$H^2_S(F_{\infty}'/F , T)$ is a torsion $R\llbracket
\Gal(F_{\infty}'/F^{\cyc})\rrbracket$-module if and only if
$H^2_S(F_{\infty}/F , T)$ is a torsion $R\ps{H}$-module and
$$\displaystyle\sum_{i\geq 1} (-1)^i\rank_{R\ps{H/N}}H_i\big(N,
H^2_S(F_{\infty}/F, T)\big) =0.$$ But this is immediate from Theorem
\ref{torsion main theorem/lemma}. \epf

We end the section with the following remark.

\begin{remark}
The results here can be extended to
the case when the ring $R$ is a commutative Noetherian local domain
which is finite flat over a regular local ring $R_0$. In this case,
for a compact pro-$p$ $p$-adic Lie group without $p$-torsion, the
notion of torsion and pseudo-nullity over $R\ps{G}$ is defined
similarly as in the regular case. It then follows from the flatness
condition of $R$ that one has the following isomorphisms
\[ R\ot_{R_0}\Ext^{i}_{R_0\ps{G}}(M , R_0\ps{G}) \cong \Ext^{i}_{R\ps{G}}(R\ot_{R_0}M ,
R\ps{G}) \cong \Ext^{i}_{R\ps{G}}(M^d , R\ps{G}), \]
 where $d$ is
the $R_0$-rank of $R$, and a $R$-free module $T$ is clearly still
$R_0$-free. Therefore, the question of pseudo-nullity over $R\ps{G}$
is reduced to pseudo-nullity over $R_0\ps{G}$, and the results in
this section apply.
\end{remark}

\section{Some examples} \label{example}

In this section, we will discuss some numerical examples of
pseudo-nullity.

(a) The first example we consider is the case $T=\Zp(1)$.  Take $F
=\Q(\mu_{p})$. Let $S$ be the set of prime(s) of $F$ above $p$. By
the theorem of Ferrero-Washington, the group
$\Gal(K(F^{\cyc})/F^{\cyc})$ is finitely generated over $\Z_{p}$.
Now if $p<1000$, it follows from \cite[Theorem 1.4]{Sh2} that
$\Gal(K(F_{\infty})/F_{\infty})$ is a finitely generated
$\Zp\ps{\Gal(F_{\infty}/F_{\infty}')}$-module (and hence a
pseudo-null $\Zp\ps{\Gal(F_{\infty}/F)}$-module) for every
$S$-admissible $p$-adic extension $F_{\infty}$ of $F$ which contains
$F_{\infty}':=\Q(\mu_{p^{\infty}}, p^{-p^{\infty}})$. Note that
$\Q(\mu_{p^{\infty}}, p^{-p^{\infty}})$ (and hence any $F_{\infty}$
containing it) satisfies $(\mathbf{Dim}_{S})$ by \cite[Lemma
3.9]{HV}. Write $N=\Gal(F_{\infty}/F_{\infty}')$. It follows that
$\Gal(K(F_{\infty})/F_{\infty})$ is finitely generated over
$\Zp\ps{N}$ if and only if $H^2_S(F_{\infty}/F, \Zp(1))$ is finitely
generated over $\Zp\ps{N}$. Therefore in this case, we have
$$ \rank_{\Zp\ps{H/N}}H_i\big(N, H^2_S(F_{\infty}/F, \Zp(1))\big)
=0$$ for all $i$ and the equation of Lemma \ref{relative rank} is
vacuous here. We mention in passing that if $p$ is a regular prime,
it will follow from an application of a classical result of Iwasawa
that $\Gal(K(F_{\infty})/F_{\infty}) =0$ for every $S$-admissible
$p$-adic extension $F_{\infty}$ of $F$ (for instances, see
\cite[Section 4]{Oc}).

(b) The next example comes from \cite[Example 23]{Bh} which is also
considered in \cite[Example 4.8]{CS}. Let $E$ be the elliptic curve
$150A1$ of Cremona's table which is given by
\[ y^2 + xy = x^3 - 3x -3,\]
Take $p=5$ and $F = \Q(\mu_5)$. Let $S$ be the set of primes of $F$
lying above $2, 3, 5$ and $\infty$. The elliptic curve $E$ has good
ordinary reduction at the unique prime of $F$ above $5$ and split
multiplicative reduction at the unique primes of $F$ above 2 and 3.
It was shown in \cite[Example 23]{Bh} that $Y(T_5 E/F_{\infty})$ is
a pseudo-null $\Z_5\ps{\Gal(F_{\infty}/F)}$-module for the
$S$-admissible $5$-adic extension $F_{\infty} = \Q(E[5^{\infty}],
3^{5^{-\infty}})$. Since $\Q(E[5^{\infty}])$ satisfies
$(\mathbf{Dim}_{S})$ (cf.\ \cite[Lemma 2.8]{C}), so does
$\Q(E[5^{\infty}], 3^{5^{-\infty}})$. Applying Theorem
\ref{pseudo-null descent}, we have that $Y(T_5 E/\mathcal{L})$ is a
pseudo-null $\Z_5\ps{\Gal(\mathcal{L}/F)}$-module when $\mathcal{L}$
is one of the following $S$-admissible $5$-adic extensions:
\[\Q(E[5^{\infty}], 2^{5^{-\infty}}, 3^{5^{-\infty}}), \quad\Q(E[5^{\infty}], 3^{5^{-\infty}}, 5^{5^{-\infty}}),\quad
\Q(E[5^{\infty}], 2^{5^{-\infty}}, 3^{5^{-\infty}},
5^{5^{-\infty}})),\] \[L_{\infty}(E[5^{\infty}], 2^{5^{-\infty}},
3^{5^{-\infty}}), \quad L_{\infty}(E[5^{\infty}], 3^{5^{-\infty}},
5^{5^{-\infty}}),\quad L_{\infty}(E[5^{\infty}], 2^{5^{-\infty}},
3^{5^{-\infty}}, 5^{5^{-\infty}}), \] where $L_{\infty}$ is any
$\Z_5^r$-extension of $F$ for $1\leq r\leq 3$.

Write $K_{\infty}= \Q(E[5^{\infty}])$. Applying Theorem
\ref{pseudo-null descent} in this direction, we have
$$\rank_{\Z_5\ps{\Gal(K_{\infty}/F^{\cyc})}} Y(T_5E/K_{\infty}) =
\rank_{\Z_5\ps{\Gal(K_{\infty}/F^{\cyc})}}
H_1\big(\Gal(F_{\infty}/K_{\infty}), H^2_S(F_{\infty}/F,
T_5E)\big).$$
  Unfortunately, we do not know how to show that this latter quantity is
  zero which will then verify the pseudo-nullity for
  $Y(T_5E/L_{\infty})$. What we do know at present is that this
  quantity is either zero or two (cf. \cite[Example 4.8]{CS}).

(c) We now discuss our final numerical example which has also been
considered in \cite[p. 362]{Jh}. Let $E$ be the elliptic curve
$79A1$ of Cremona's table which is given by
\[ y^2 + xy + y = x^3 + x^2 -2x.\]
Take $p=3$ and $F = \Q(\mu_3)$. Let $S$ be the set of primes of $F$
lying above $3, 79$ and $\infty$. The elliptic curve $E$ has good
ordinary reduction at the unique prime of $F$ above $3$ and non
split multiplicative reduction at the two primes of $F$ above 79. It
was shown in \cite[p. 362]{Jh} that $Y(T_3E/F_{\infty})$ is a
pseudo-null $\Z_3\ps{\Gal(F_{\infty}/F)}$-module when $F =
\Q(\mu_{3^{\infty}}, 79^{-3^{\infty}})$. Again noting that
$\Q(\mu_{3^{\infty}}, 79^{-3^{\infty}})$ satisfies
$(\mathbf{Dim}_S)$ (cf.\ \cite[Lemma 3.9]{HV}), one can apply
Theorem \ref{pseudo-null descent} to conclude that
$Y(T_3E/\mathcal{L})$ is a pseudo-null
$\Z_3\ps{\Gal(\mathcal{L}/F)}$-module when $\mathcal{L}$ is one of
the following $S$-admissible $3$-adic extensions:
\[ \Q(\mu_{3^{\infty}},
3^{-3^{\infty}}, 79^{-3^{\infty}}), \quad
L_{\infty}(79^{-3^{\infty}}),\quad L_{\infty}(3^{-3^{\infty}},
79^{-3^{\infty}}). \] Here $L_{\infty}$ is the unique
$\Z_3^2$-extension of $F$.

Now the residual representation on the Tate module of the elliptic
curve $79A1$ is irreducible for $p=3$. Hence there exists a complete
Noetherian local domain $R$ which is finite flat over $\La =
\Z_3\llbracket X\rrbracket$ and a free $R$-module $T$ of rank $2$
with a continuous $\Gal(\bar{\Q}/\Q)$-action which is the $R$-dual
of the Galois representation attached to the Hida family associated
to the weight 2 newform corresponding to $E$ (cf.\ \cite{Hi}, also
see \cite[Section 1 and Remarks at the end of Section 6]{SS} for
details on $R$ and $T$; note that our $T$ here is the $R(1)$-dual of
the representation considered there). As mentioned above, the
$S$-admissible $3$-adic extension $\Q(\mu_{3^{\infty}},
79^{-3^{\infty}})$ satisfies $(\mathbf{Dim}_{S})$, and therefore any
$S$-admissible $3$-adic extension containing $\Q(\mu_{3^{\infty}},
79^{-3^{\infty}})$ will also satisfy $(\mathbf{Dim}_{S})$. Applying
Theorem \ref{pseudo-null main} (and noting the remarks made at the
end of Section \ref{pseudo-nullity section}), we have that
$Y(T/\mathcal{L})$ is a pseudo-null
$R\ps{\Gal(\mathcal{L}/F)}$-module when $\mathcal{L}$ is one of the
following $S$-admissible $3$-adic extensions:
\[\Q(\mu_{3^{\infty}}, 79^{-3^{\infty}}), \quad\Q(\mu_{3^{\infty}},
3^{-3^{\infty}}, 79^{-3^{\infty}}),\quad
L_{\infty}(79^{-3^{\infty}}), \quad L_{\infty}(3^{-3^{\infty}},
79^{-3^{\infty}}). \] (Alternatively, one can deduce this by
combining Jha's observation and Theorem \ref{pseudo-null descent}.)

\section{Complement: On torsionness of fine Selmer group}
\label{Torsion section}

 As before, $p$ will denote a prime. Let $F$ be a number
field, where we assume that it has no real primes when $p=2$.
 We also assume that our
Galois module $T$ is a free $R$-module of finite rank with a
continuous $R$-linear $G_S(F)$-action, where $R$ is a complete
regular local ring with finite residue field of characteristic $p$.
The following is a weaker version of the conjecture made in Section
\ref{Cyclotomic Zp extension}.

\medskip \noindent \textbf{Conjecture A$'$} \textit{For
any number field $F$ and a $S$-admissible $p$-adic Lie extension
$F_{\infty}$ of $F$,\, $Y_S(T/F_{\infty})$ is a finitely generated
torsion $R\ps{G}$-module.}

\medskip
It is clear that Conjecture A$'$ will follows from Conjecture A and
Lemma \ref{fg La H}. In particular, by Theorem \ref{main},
Conjecture A$'$ is a consequence of the Iwasawa $\mu$-conjecture. We
now record the following lemma (compare with \cite[Lemma 3.1]{CS}).
We write $A= T\ot_R R^{\vee}$.

\bl \label{torsion H2 Y} Let $F_{\infty}$ be a $S$-admissible
$p$-adic Lie extension of $F$. Then the following statements are
equivalent.

\smallskip
$(a)$ $H^2_{S}(F_{\infty}/F, T)$ is a torsion $R\ps{G}$-module.

\smallskip
$(b)$ $Y_S(T/F_{\infty})$ is a torsion $R\ps{G}$-module.

\smallskip
$(c)$ $H^2(G_S(F_{\infty}),A) = 0$. \el

\begin{remark}
In view of Statement (c), Conjecture A$'$ is also sometimes called
the ``weak Leopoldt conjecture for $A$'' (over $F_{\infty}$).
\end{remark}

\bpf The equivalence of (a) and (b) follows from a similar argument
as in Lemma \ref{fg H2}. We will now establish the equivalence of
(a) and (c). Consider the following general version of the spectral
sequence of Jannsen \cite[Theorem 1]{Ja}
  $$ E_2^{i,j} = \Ext^i_{R\ps{G}}\big( H^j(G_S(F_{\infty}),A)^{\vee},
  R\ps{G}\big)\Longrightarrow H^{i+j}_S(F_{\infty}/F, T).
  $$
  (One can obtain this spectral sequence by combining the middle two derived isomorphisms in
\cite[1.6.12(4)]{FK}.) 
    Since the spectral sequence is bounded (as $R\ps{G}$ has finite projective dimension),
    it follows that $E_m^{r,s}$ must stabilize for large enough $m$. In particular, we have that $E_{\infty}^{i,j}$
    is a subquotient of $E_2^{i,j}$. By the Auslander
    regularity of $R\ps{G}$, the terms $E_2^{i,j}$, and hence $E_{\infty}^{i,j}$, are torsion $R\ps{G}$-modules for $i\neq
    0$. Since $H^{2}_S(F_{\infty}/F, T)$ has a finite filtration with
    factors $ E_{\infty}^{i,2-i}$ for $i =0,1,2$, this in turn yields
     $ \rank_{R\ps{G}}H^{2}_S(F_{\infty}/F, T) = \rank_{R\ps{G}}E_{\infty}^{0,2}.$
    On the other hand, one sees easily that the edge map
    $E_{\infty}^{0,2}\to E_2^{0,2}$ is injective and has cokernel isomorphic to a subquotient of
     \[
        \Ext_{R\ps{G}}^{2}(H^1(G,A)^{\vee}, R\ps{G})\oplus \Ext_{R\ps{G}}^{3}(H^0(G,A)^{\vee},
        R\ps{G}).
    \] As observed above, these are torsion over $R\ps{G}$. Therefore, we may
    conclude that
     $$ \rank_{R\ps{G}}H^{2}_S(F_{\infty}/F, T) = \rank_{R\ps{G}}\Hom_{R\ps{G}}\big( H^2(G_S(F_{\infty}),A)^{\vee},
  R\ps{G}\big). $$
  Therefore, it follows that $H^2_{S}(F_{\infty}/F, T)$ is a torsion $R\ps{G}$-module
  if and only if \linebreak $\Hom_{R\ps{G}}\big( H^2(G_S(F_{\infty}),A)^{\vee},
  R\ps{G}\big)$
  is a torsion $R\ps{G}$-module.
  Since $H^2(G_S(F_{\infty}),A)^{\vee}$ is reflexive by
  \cite[Proposition 3.5(ii)]{SS} and \cite[Proposition 3.11(i)]{V1},
  the latter statement holds if and only if $H^2(G_S(F_{\infty}),A) = 0$. \epf

Combining the preceding lemma with Theorem \ref{torsion main
theorem/lemma} and Corollary \ref{torsion x corollary}, we have the
following results for the torsion property of the dual fine Selmer
groups which is analogous to the pseudo-nullity results obtained in
Section \ref{pseudo-nullity section}.

\bp \label{torsion descent} Let $F_{\infty}$ be a $S$-admissible
$p$-adic Lie extension of $F$. Suppose that $F_{\infty}'$ is another
$S$-admissible $p$-adic Lie extension of $F$ which satisfies the
following properties.

\smallskip
$(i)$  $F'_{\infty}$ is contained in $F_{\infty}$.

\smallskip
$(ii)$ The group $N:=\Gal(F_{\infty}/F_{\infty}')$ is a solvable
uniform pro-$p$ group and $G/N$ has no $p$-torsion.

Then $Y_S(T/F_{\infty}')$ is a torsion $R\ps{
\Gal(F_{\infty}'/F)}$-module if and only if $Y_S(T/F_{\infty})$ is a
torsion $R\ps{G}$-module and $\displaystyle\sum_{i\geq 1}
(-1)^i\rank_{R\ps{G/N}}H_i\big(N, H^2_S(F_{\infty}/F, T)\big)
=0$.\ep

\bp \label{torsionness x} Let $F_{\infty}$ be a $S$-admissible
$p$-adic Lie extension of $F$. Let $x$ be a nonzero and nonunital
element of $R$, and suppose that $\bar{R} := R/xR$ is also regular
local. Write $\bar{T} = T/xT$. Then $Y_S(\bar{T}/F_{\infty})$ is a
torsion $\bar{R}\ps{G}$-module if and only if $Y_S(T/F_{\infty})$ is
a torsion $R\ps{G}$-module and $H^2_{S}(F_{\infty}/F, T)[x]$ is a
torsion $\bar{R}\ps{G}$-module. \ep

\begin{remark}
Proposition \ref{torsion descent} can be viewed as a somewhat
analogous statement for Selmer groups as in \cite[Theorem 2.8]{HV}
and \cite[Theorem 2.3]{HO}.
\end{remark}

We end by briefly mentioning some known cases of Conjecture A$'$. In
the case when $T=\Zp(1)$ and $F_{\infty} = F^{\cyc}$, this follows
from a classical theorem of Iwasawa \cite[Theorem 5]{Iw}. For a
general admissible extension $F_{\infty}$, one can deduce the
conjecture from the cyclotomic case using a limit argument with
statement (c) of Lemma \ref{torsion H2 Y} (see \cite[Theorem
6.1]{OcV} for details).

When $T$ is the Tate module of an elliptic curve $E$ which does not
have potentially supersingular reduction at any primes above $p$, it
is in fact conjectured that the Pontryagin dual of the classical
Selmer group is a torsion $\Zp\ps{G}$-module (see \cite{HO, Maz,
Sch}). Conjecture A$'$ is then a consequence of this more general
conjecture.

In the case when $E$ is an elliptic curve defined over $\Q$ with
good ordinary reduction at $p$ and $F$ is an abelian extension of
$\Q$, it follows from a deep theorem of Kato \cite{K} that the dual
Selmer group over $F^{\cyc}$ is a torsion $\Zp\ps{\Ga}$-module. In
particular, conjecture A$'$ holds in this case. Note that one cannot
apply a limit argument as in \cite[Theorem 6.1]{OcV} in this
situation to the case of a general admissible $p$-adic Lie extension
$F_{\infty}$, since in general a finite extension $L$ of $F$
contained in $F_{\infty}$ need not be abelian over $\Q$ and
therefore, Kato's theorem does not apply. However, one may still
apply Proposition \ref{torsion descent} to obtain cases of
Conjecture A$'$ over solvable admissible $p$-adic Lie extension of
$F$.

On the other hand, one can deduce conjecture A$'$ for elliptic curve
over a $S$-admissible $p$-adic Lie extension of $F$ containing
$F(E[p^{\infty}])$ via \cite[Lemma 2.4]{CS}. Here one does not need
any reduction hypothesis on $E$ nor assumption on $F$.

Finally, we observe that under appropriate modification as noted in
the Remark at the end of Section \ref{pseudo-nullity section}, we
can extend \cite[Corollary 6.17]{SS} to a solvable admissible
$p$-adic Lie extension.

\appendix
\renewcommand\thesection{Appendix \Alph{section}}
\section{On the structure of $R\ps{G}$}
\renewcommand\thesection{\Alph{section}}

The purpose of this appendix is to prove certain results on the
structure of the completed group algebra $R\ps{G}$. We believe
results of such are well-known among experts, although they do not
seem to have been written down properly anywhere except for the case
$R=\Zp$. Since these results are basic for the discussion in this
paper, we have included their proofs here. The following is the main
result of the appendix.

\bt \label{compact Auslander regular} $(a)$ If $R$ is a commutative
Noetherian local ring with finite residue field of characteristic
$p$ and $G$ is a compact $p$-adic Lie group, then $R\ps{G}$ is
Noetherian.

$(b)$ If $R$ is a commutative Noetherian local domain with finite
residue field of characteristic $p$ and $G$ is a uniform pro-$p$
group, then $R\ps{G}$ has no zero divisor.

$(c)$ If $R$ be a commutative complete regular local ring with
finite residue field of characteristic $p$ and $G$ is a compact
$p$-adic Lie group without $p$-torsion, then $R\ps{G}$ is an
Auslander regular ring. \et

Statement (a) and (b) are well-known theorems of Lazard \cite{Laz}
when $R=\Zp$ or $R$ is a finite field of order $p$. (see also
\cite[Corollary 7.25, Corollary 7.26]{DSMS}). Statement (a) in this
generality has been established in \cite[Theorem 8.7.8]{Wil} and
\cite[Proposition 3.0.1]{LS}. Neumann \cite{Neu} has shown statement
(b) when $G$ is a pro-$p$ $p$-adic Lie group without $p$-torsion and
$R=\Zp$ by a different approach. His method can be easily extended
to the case when $R$ is a regular local ring with characteristic
zero. However, his method does not seem to apply to the case when
$R$ has characteristic $p$. Statement (c) is the theorem of Venjakob
when $R=\Zp$ or $R$ is a finite field of order $p$ \cite[Theorems
3.26, 3.30(b)]{V1}. In this appendix, we will give a
uniform\footnote{literally and figuratively} approach to prove all
three statements simultaneously.

We first review some facts which can be found in \cite{DSMS}. For a
finitely generated pro-$p$ group $G$, we write $G^{p^i} = \langle
g^{p^i}|~g\in G\rangle$, that is, the group generated by the
$p^i$th-powers of elements in $G$. The pro-$p$ group $G$ is said to
be \textit{powerful} if $G/\overline{G^{p}}$ is abelian for odd $p$,
or if $G/\overline{G^{4}}$ is abelian for $p=2$. We define the lower
$p$-series by $P_{1}(G) = G$, and
\[ P_{i+1}(G) = \overline{P_{i}(G)^{p}[P_{i}(G),G]}, ~\mbox{for}~ i\geq 1. \] It follows from \cite[Thm.\
3.6]{DSMS} that if $G$ is a powerful pro-$p$ group, then $G^{p^i} =
P_{i+1}(G)$. This in turn implies that a finitely generated pro-$p$
group $G$ is powerful if and only if $[G,G]\sbs G^{p^2}$. It follows
from the same theorem that for a finitely generated powerful pro-$p$
group $G$, the $p$-power map
\[ P_{i}(G)/P_{i+1}(G)\stackrel{\cdot p}{\lra}
P_{i+1}(G)/P_{i+2}(G)\] is surjective for each $i\geq 1$. If the
$p$-power maps are isomorphisms for all $i\geq 1$, we say that $G$
is \textit{uniformly powerful} (abrev.\ \textit{uniform}). Note that
in this case, we have an equality $|G:P_2(G)| = |P_i(G):
P_{i+1}(G)|$ for every $i\geq 1$. We now recall the following
characterization of compact $p$-adic Lie groups due to Lazard
\cite{Laz} (see also \cite[Cor.\ 8.34]{DSMS}): a topological group
$G$ is a compact p-adic Lie group if and only if $G$ contains a open
normal uniform pro-$p$ subgroup.

For the remainder of the appendix, $G$ will always be a uniform
pro-$p$ group, unless otherwise stated. For $n\geq 0$, we shall
write $G_n = P_n(G)$. Denote $I_n$ to be the augmentation kernel of
the map $R\ps{G}\twoheadrightarrow R[G/G_n]$ which is a closed
two-sided ideal of $R\ps{G}$. Since the subgroups $G_n$ form a basis
of neighbourhood of $1$ in $G$, we have $R\ps{G} = \plim_n
R\ps{G}/I_n$. Fix a minimal set of topological generators $a_{1},
a_{2},...,a_{d}$ of $G$. Write $b_{i}=a_{i} -1$ for each $i$. For
$\al = (\al_{1},...,\al_{d})\in\N^{d}$ (here, $\N$ is the set of
natural numbers including $0$) and any $d$-tuple $\mathbf{v} =
(v_{1},...,v_{d})\in \La^d$, we write
\[ \langle \al\rangle = \al_{1} + \cdots + \al_{d},\quad\mathbf{v}^{\al} =
v_{1}^{\al_{1}}\cdots v_{d}^{\al_{d}}. \] In particular, we write
$\bb^{\al} = b_{1}^{\al_{1}}\cdots b_{d}^{\al_{d}}$. We now examine
the structure of $R\ps{G}$ as an $R$-module. For $n\geq 1$, we
define
\[ T_{n} =\{\al\in\N^{d} ~|~
\al_{i}<p^{n-1}~\mathrm{for}~i=1,...,d\}.\]

\bl \label{sum decomposition of group algebra}  If $G$ is uniform,
then we have a direct sum decomposition \[ R\ps{G} = I_{n}
\oplus\bigoplus_{\al\in T_{n}}R \bb^{\al}
\] of $R$-modules. Furthermore, we have $R\bb^{\al}\cong R$ for each $\al$ and
\[ I_n = I_{n+1} \oplus\bigoplus_{\al\in T_{n+1}\setminus T_{n}}R \bb^{\al} \] \el

\bpf (Compare with \cite[Lemma 7.9]{DSMS}) Denote $\phi$ to be the
canonical quotient $R\ps{G}\tha R[G/G_{n}]$. By \cite[Theorem
3.6]{DSMS}, every element of $G/G_{n}$ can written as
$a_{1}^{\al_{1}}\cdots a_{d}^{\al_{d}}G_{n}$ with $\al_{i}<
p^{n-1}$. Hence $\{\phi(\mathbf{a}^{\al})|~\al\in T_{n}\}$ generates
$R[G/G_{n}]$ as an $R$-module. Since $G$ is uniform, we have
$|G/G_{n}|= p^{(n-1)d}$, and so $\phi(R\ps{G})$ is a free $R$-module
of rank $p^{(n-1)d}$. On the other hand, we have $|T_{n}| =
p^{(n-1)d}$. Therefore, the generating set
$\{\phi(\bb^{\al})|~\al\in T_{n}\}$ is actually a free $R$-basis for
this module. The assertions in the lemma are now immediate from
this. \epf

The next lemma is the key ingredient to the proof of our theorem. We
write $\mathrm{gr}_{\{I_n\}}R\ps{G}$ for the graded ring
$\bigoplus_{n\geq 0}I_n/I_{n+1}$ whose multiplication is given by
 $$(\mathbf{b}^{\al} + I_{n+1})(\mathbf{b}^{\be} + I_{m+1}) =
\mathbf{b}^{\ga} + I_{n+m+1},$$
 where $\ga_i = \al_i + \be_i$. To see that this multiplication is well-defined,
 we note that by an induction argument, it suffices to show that $b_i b_j - b_j b_i
 \in I_2$. Now observe that
\[ b_i b_j - b_j b_i = a_i a_j - a_j a_i = ([a_i, a_j]-1)a_j a_i\]
Since $G$ is powerful, we have that $[a_i, a_j] \in G_2$. This
implies that $[a_i, a_j]-1 \in I_2$ which in turn implies that $b_i
b_j - b_j b_i \in I_2$, as required. We may now prove our lemma.

\bl \label{uniform bijective} Let $R$ be a complete Noetherian local
ring with finite residue field of characteristic $p$. Let $G$ be a
uniform pro-$p$ group. Then we have an
 $R$-algebra isomorphism given by
\[ \ba{c} \Phi : R[X_{1},...,X_{d}]\lra \mathrm{gr}_{\{I_n\}}R\ps{G}\\
 X_j\mapsto b_{j} + I_2 . \ea\] \el

\bpf The above argument shows that the assignment $X_{j}\mapsto
b_{j} + I_2$ give a well-defined $R$-algebra homomorphism.
Surjectivity of $\Phi$ is an immediate consequence of Lemma \ref{sum
decomposition of group algebra}. It remains to show that $\Phi$ is
injective. Let $A_n$ be the $R$-submodule of $R[X_{1},...,X_{d}]$
generated by monomials of the form $X_1^{\al_1}\cdots X_d^{\al_d}$,
where $\al_i <p^{n-1}$ for $i =1,...,d$. Then one has
\[ \Phi( A_n) = \bigoplus_{j=1}^n I_j/I_{j+1} = \bigoplus_{\al\in
T_n}R\bb^{\al}, \] where the second equality comes from Lemma
\ref{sum decomposition of group algebra}. Now $A_n$ and $\Phi(A_n)$
are both free $R$-modules of rank $m:= |T_n| = p^{(n-1)d}$. Choose
$R$-isomorphisms $f : R^m \stackrel{\sim}{\lra} A_n$ and $g : R^m
\stackrel{\sim}{\lra} \Phi(A_n)$. Then $g^{-1}\circ\Phi\circ f$ is a
surjective endomorphism of $R^m$. By \cite[Theorem 2.4]{M}, it
follows that $g^{-1}\circ\Phi\circ f$ is an automorphism. In
particular, this implies that $\ker\Phi|_{A_n} =0$. Therefore, we
have $\ker\Phi = \cup_n \ker\Phi|_{A_n} =0$, hence proving the
lemma. \epf

We can now give a proof of Theorem \ref{compact Auslander regular}.

\bpf[Proof of Theorem \ref{compact Auslander regular}] Statement (b)
is an immediate consequence of Lemma \ref{uniform bijective} and
\cite[Proposition 7.27(i)]{DSMS}. Now let $U$ be an open normal
uniform pro-$p$ subgroup of $G$. Since $U$ is a subgroup of $G$ with
finite index, it follows that $R\ps{G}$ is Noetherian if $R\ps{U}$
is so. The latter is then an immediate consequence of Lemma
\ref{uniform bijective} and \cite[Proposition 7.27(ii)]{DSMS}.

It remains to show statement (c). By (a), the ring $R\ps{G}$ is
Noetherian. Therefore, the global dimension of $R\ps{G}$ coincides
with its topological projective dimension. Let $U$ be an open normal
uniform pro-$p$ subgroup of $G$.  Since $G$ has no $p$-torsion, we
may apply a result of Serre (cf.\ \cite[Corollaire 1]{Se}) to
conclude that $\cd_p (G) =\cd_p (U)$. It then follows from
\cite[Theorem 4.1]{Bru} that $R\ps{G}$ has finite global dimension.
It remains to verify the Auslander condition. By \cite[Proposition
5.4.17]{NSW}, we have an isomorphism $\Ext^{i}_{R\ps{G}}(M, R\ps{G})
\cong \Ext^{i}_{R\ps{U}}(M, R\ps{U})$ of $R\ps{U}$-modules for any
$R\ps{G}$-module $M$. Therefore, we are reduced to showing that
$R\ps{U}$ is Auslander regular. By Lemma \ref{uniform bijective}, we
have that the associated graded ring of
$\mathrm{gr}_{\{I_n\}}R\ps{U}$ is a commutative regular local ring,
since it is a polynomial ring over a regular local ring. The
conclusion will follow from an application of the theorem of
Bj$\ddot{\mathrm{o}}$rk (see Remarks after \cite[Theorem 3.9]{Bjo}
or \cite[Theorem 3.21]{V1}). (One still need to verify the closure
condition in the cited theorem but this follows immediately from the
observation that the $I_n$'s are two-sided closed ideals of
$R\ps{U}$). \epf

 \small


\end{document}